\documentclass[a4paper,11pt]{article}

\usepackage{amsmath}
\usepackage{amsfonts}
\usepackage{amssymb}
\usepackage{amsthm}
\usepackage{graphicx}
\usepackage{verbatim}
\usepackage{authblk}
\usepackage{setspace} 

\newcommand{\eps}{\varepsilon}
\newcommand{\N}{\mathbb{N}}
\newcommand{\Q}{\mathbb{Q}}
\newcommand{\R}{\mathbb{R}}
\newcommand{\T}{\mathbb{T}}
\newcommand{\Z}{\mathbb{Z}}
\newcommand{\RR}{\mathcal{R}}
\newcommand{\TT}{\mathcal{T}}
\newcommand{\deltabf}{\text{\boldmath $\delta$}}

\date{}
\author[1]{Pau Rabassa}
\author[2]{Angel Jorba}
\author[2]{Joan Carles Tatjer}
\affil[1] {\mbox{Johann Bernoulli Institute for Mathematics and Computer Science,}
\centerline{ \mbox{University of Groningen, Groningen, The Netherlands}}
\newline 
\mbox{E-mail: {\tt paurabassa@gmail.com}} 
\vspace{2mm}}

\affil[2]{\mbox{Departament of Matem\`atica Aplicada i An\`alisi,}
\mbox{Universitat de Barcelona, Barcelona, Spain} 
\mbox{E-mails: {\tt angel@maia.ub.edu}, {\tt jcarles@maia.ub.es}}} 

\title{
{N}umerical evidences of universality and self-similarity 
in the Forced Logistic Map\thanks{This work has been supported 
   by the MEC grant MTM2009-09723 and the CIRIT grant 2009 SGR 67. 
   P.R. has been partially supported by the PREDEX project, funded 
   by the Complexity-NET: {\tt www.complexitynet.eu}.}
}

\begin{document}

\maketitle
\begin{abstract}
We explore  different families of quasi-periodically Forced
Logistic Maps for the existence of universality and 
self-similarity properties. In the bifurcation diagram 
of the Logistic Map it is well known that there exist 
parameter values $s_n$ where the $2^n$-periodic orbit is 
superattracting. Moreover these parameter values lay between one 
period doubling and the next. 
Under quasi-periodic forcing, the superattracting periodic orbits 
give birth to two reducibility-loss 
bifurcations in the two dimensional parameter space of the Forced 
Logistic Map, both around the points $s_n$. 
In the present work we study numerically the asymptotic 
behavior of the slopes of these bifurcations with respect to $n$. 
This study evidences the existence of universality properties 
and self-similarity of the bifurcation diagram in the parameter space.
\end{abstract}


{\bf Universality and self-similarity properties of uniparametric 
families of unimodal maps are a well known phenomenon. The paradigmatic 
example of this phenomenon is the Logistic Map $l_\alpha(x) = \alpha x(1-x)$. 
Given a typical one parametric family 
of unimodal maps $\{l_\alpha\}_{\alpha\in I}$  one  observes
numerically that there exists a sequence of parameter
values $\{d_n\}_{n\in \N} \subset I$ such that
the attracting periodic orbit of the map undergoes a period doubling
bifurcation. Between one period doubling and the next there 
exists a parameter value $s_n$, for which the critical point
of $l_{s_n}$ is a periodic orbit with period $2^n$. One can
also observe that
\begin{equation}
\label{universal limit sumicon}
\lim_{n\rightarrow \infty} \frac{d_n - d_{n-1}} {d_{n+1} - d_{n}} = 
\lim_{n\rightarrow \infty} \frac{s_n - s_{n-1}} {s_{n+1} - s_{n}} =
\deltabf  = \texttt{ 4.66920...}. 
\end{equation}
The convergence to this limit $\deltabf$ (the 
so-called Feigenbaum constant) indicates a self-similarity 
on the parameter space of the family. On the other hand, 
the constant $\deltabf$ is universal,
in the sense that  
one obtains the same ratio $\deltabf$
for any family  of unimodal maps with a quadratic
turning point having a cascade of period doubling bifurcations.
In this paper we explore if the same kind of phenomenon can be 
observed when the  one dimensional maps is forced 
quasi-periodically. The answer is that universality and 
self-similarity do manifest, but they do it in different 
manners. Moreover, we show that they occur in a more
``restrictive'' class of maps, in the sense that the 
quasi-periodic forcing has to have a very particular 
form. We provide a theoretical explanation to this 
phenomenon in \cite{JRT11a, JRT11b, JRT11c} 
(see also \cite{Rab10}). }

\section{Introduction}
\label{sec:Intro}

In the late seventies, Feigenbaum (\cite{Fei78,Fei79}) and 
 Coullet and Treser (\cite{CT78}) proposed at the same 
time the renormalization operator to explain the universal 
features observed in the cascade of
period doubling bifurcations of the Logistic Map. 
This explanation was based on the existence of a 
hyperbolic fixed point of the operator with suitable 
properties. The first proof of the existence of this 
point and his hyperbolicity were obtained with numeric
assistance \cite{Lan82,EW87}. In \cite{Sul92} Sullivan
generalized the operator and provided a theoretical 
proof of the hyperbolicity using complex dynamics. In \cite{Lyu99, 
MvS93} a nice summary of the one dimensional 
renormalization theory can be found, up to the date of their respective
publication, as well as contributions to this theory.

The kind of maps that we consider are maps 
in the cylinder where the dynamics on the periodic variable 
is given by a rigid rotation and the dynamics on the other variable 
is given by an unimodal one dimensional map plus a small perturbation 
which depends on both variables.  These maps are known as 
quasi-periodically forced one dimensional maps and they have been 
extensively studied (\cite{GOPY84,HH94, PMR98,Jag03,FKP06,JT08,Bje09,FH10})
with a focus on the existence of strange non-chaotic attractors. 

The paradigmatic example in our case of study is the Forced 
Logistic Map, which is the map on the 
cylinder $\T\times \R$ defined as
\begin{equation}
\label{FLM}
\left.
\begin{array}{rcl}
\bar{\theta} & = & \theta + \omega  ,\\
\bar{x} & = & \alpha x(1-x)(1+ \eps \cos(2\pi \theta)),  \end{array}
\right\}
\end{equation}
where $(\alpha,\eps)$ are parameters and $\omega$ is a fixed Diophantine
number (typically it will be the golden mean). This 
family has interest not only for phenomena related with 
the existence of strange non-chaotic attractors 
(\cite{HH94,PMR98,FKP06}) but also 
as a toy-model for the truncation of period doubling bifurcation 
cascade (\cite{Kan83,JRT11p}). The term ``truncation of the period 
doubling bifurcation cascade'' refers to the fact that, 
when one fixes $\eps$ and let $\alpha$ grow, the attracting 
set of the map undergoes only a finite number of 
period doubling bifurcations before exhibiting a chaotic behavior. 
This differs from the one dimensional case, where the number of period 
doubling bifurcations before chaos is infinite. 

\begin{figure}[t]
\begin{center}
\includegraphics[width=12cm]{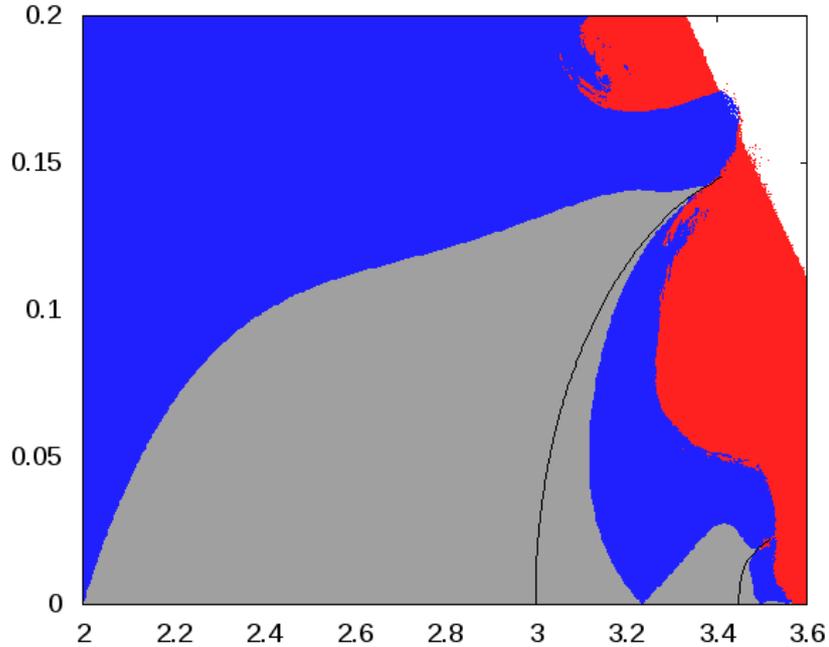}
\caption{ 
Diagram of the parameter space of the map (\ref{FLM}) for
$\omega= \frac{\sqrt{5}-1}{2}$. The axis correspond to
the parameters $\alpha$ (horizontal) and $\eps$ (vertical).
For the correspondence of each color with
the properties of the attractor see Table \ref{table coding diagrams}.
The black lines represent the period doubling bifurcations from 
period one to two (left) and two to four (right).}
\label{FLM parameter space}
\end{center}
\end{figure}

In \cite{JRT11p} we studied the truncation of the period 
doubling cascade for the map (\ref{FLM}). We observed that the reducibility 
of the attracting set plays a crucial role. We computed 
bifurcation diagrams in terms
of the dynamics of the attracting set, taking
into account different properties,  as the Lyapunov exponent
and, in the case of having a periodic invariant curve, its
period and its reducibility. One of these bifurcation diagrams 
is reproduced in Figure \ref{FLM parameter space} (see Table 
\ref{table coding diagrams} for the label of every color).

\begin{table}[t]
\begin{center}
\begin{tabular}{|lcl|}
\hline
\rule{0pt}{3ex} Color  & \phantom{c} & Dynamics of the attractor \\ \hline
\rule{0pt}{3ex} Black  &\phantom{c} & Invariant curve with zero Lyapunov exponent \\
\rule{0pt}{2ex} Red &\phantom{c} & Chaotic attractor  \\
\rule{0pt}{2ex} Blue & \phantom{c} & Non-chaotic non-reducible attractor \\
\rule{0pt}{2ex} Grey &\phantom{c} & Non-chaotic reducible attractor \\
\rule{0pt}{2ex} White  &  \phantom{c} & No attractor (divergence to $-\infty$) \\
\hline
\end{tabular}
\end{center}
\vspace{-4mm}
\caption{Color coding for the Figure
\ref{FLM parameter space}. }
\label{table coding diagrams}
\end{table}

Let $d_n$ be the parameter value where the attracting periodic 
orbit of the one dimensional map doubles from period $2^n$ to period 
$2^{n+1}$. Figure \ref{table coding diagrams} reveals that from every 
parameter $(\alpha,\eps) = (d_n,0)$ of the map (\ref{FLM}) it is born a period 
doubling bifurcation curve of the attracting set. Let $s_n$ be the 
parameter value where the critical point of the (non-forced) one 
dimensional family is periodic with period $2^n$.
In the case of analytic maps in the cylinder, the reducibility 
loss of an invariant curve can be characterized as a bifurcation, 
for more details see definition 2.3 in \cite{JRT11p}. 
Figure \ref{table coding diagrams} also reveals that 
from every parameter value $(\alpha,\eps) = (s_n, 0)$
two curves of reducibility loss are born. These curves correspond to a
reducibility-loss and ``reducibility-recovery'' bifurcations 
of the $2^n$-periodic invariant curve. In \cite{JRT11a} we 
prove, under suitable hypothesis, that these curves exist. 
In Figure \ref{table coding diagrams} we can observe that the 
period doubling bifurcation curve born around $(d_n,0)$ is confined 
by one of the reducibility loss bifurcation curves born around 
$(s_n,0)$ and  another one born around $(s_{n+1},0)$. 

In this paper we approximate numerically slopes of the reducibility 
loss bifurcations curves described above and use them to explore for the 
existence of universality and self-similarity 
properties within the quasi-periodic one dimensional maps. 
A posteriori we know that the self-similarity properties exists 
between the bifurcation diagrams of the map (\ref{FLM}) for different
values of $\omega$. For a better presentation of the concept, let us
start the discussion looking for self-similarity withing the bifurcation diagram 
of the map (\ref{FLM}) with $\omega$ fixed.

\section{Description of the computations} 
\label{Sec:Descrip Comp}

Consider the sequence of values $s_n$ for which the Logistic Map has a
superstable periodic orbit as before. In Figure \ref{FLM parameter space} 
we can observe how the two curves (namely $S_n^-$ and $S_n^+$) that 
are born at the parameter values $(\alpha,\eps)= (s_n, 0)$ define 
a region of non-reducibility around this point. Assume that each of these
curves correspond to reducibility-loss bifurcations and 
they can be written locally with $\alpha$ as a graph of
$\eps$ (in \cite{JRT11a} we give concrete conditions under 
which this is true). 
In other words, there exist a neighborhood $U_n$ of $s_n$, 
interval $I=[0,\rho)$ and functions $\alpha: I\rightarrow \R$ 
and $\beta: I\rightarrow \R$ such that $S^-_n\cap U= 
\{(\alpha_n(\eps),\eps) | \thinspace \eps\in I\}$ 
and $S^+_n\cap U= \{(\beta_n(\eps),\eps) | \thinspace \eps\in I\}$.

Assume that the functions $\alpha_n(\eps)$ and $\beta_n(\eps)$ are 
written at first order as
\[\alpha_n(\eps) = s_n + \alpha'_n \eps + o(\eps), 
\text{ and } \beta_n(\eps) = s_n + \beta'_n \eps + o(\eps).
\]

Let $s^*$ be the limit of the parameter values $s_n$. The self-similarity 
in the bifurcation diagram of the Logistic Map is manifested in the 
following way. One observes that $s_n -s^* \approx \delta(s_{n+1} - s^*)$. 
In terms of the parameter
space this corresponds to the fact that the
affine map $L(\alpha) = \delta(\alpha- s^*) + s^*$ sends (approximately)
the point $s_{n}$ to $s_{n-1}$. We would like to find an analog
self-similarity in the parameter space of the Forced Logistic Map 
(\ref{FLM}). We can use the curves $S^\pm_n$ to detect this self-similarity. 
Therefore we look for an affine map of the kind
\begin{equation}
\label{equation afine transformation parameters}
L\left(\begin{array}{c} \alpha \\ \eps \end{array} \right) = 
\left(\begin{array}{cc} \delta_0 & 0 \\ 0 & \delta_1 \end{array} \right) 
\left(\begin{array}{c} \alpha - s^*  \\ \eps \end{array} \right)  
+ \left(\begin{array}{c} s^* \\ 0  \end{array} \right),  
\end{equation}
such that it maps the curves $S^-_{n}$ to $S^-_{n-1}$ (and respectively
$S^+_{n}$ to $S^+_{n-1}$). If we impose these conditions to the local
parameterization (around $s_i$) of the curves considered above
then we get
\[
\left( \begin{array}{c}  \alpha_{n-1} + \alpha_{n-1}' t + o(t) \\ t \end{array} \right) = 
\left(\begin{array}{cc} \delta_0 & 0 \\ 0 & \delta_1 \end{array} \right)
\left( \begin{array}{c}  \alpha_n + \alpha_n' s + o(s) - s^* \\ s \end{array} \right) 
+ \left(\begin{array}{c} s^* \\ 0  \end{array} \right) , 
\]
for any $t$ and $s$. Then, replacing $s=t/\delta_1$ in the first coordinate
and equating terms in the order of $t$, we obtain
\[
\left\{ 
\begin{array}{rcl} \alpha_n & =&  \delta_0 \alpha_{n-1} + (1-\delta_0) s^* \\ 
\alpha'_n & = & \frac{\delta_0}{\delta_1} \alpha'_{n+1}. \end{array} 
\right. 
\]

Using these two equations we have that the value $\delta_1$ can be
estimated as $\delta_1 \approx \delta_0 \frac{\alpha'_{n}}{\alpha'_{n-1}}$.
For $\delta_0$ we recover the estimation 
$\delta_0 \approx \frac{\alpha_{n} - \alpha_{n-1}}{\alpha_{n-1} - \alpha_{n-2}}$,
which converges to the Feigenbaum constant $\deltabf$. Therefore we 
replace $\delta_0$ by $\deltabf$ in the estimation of $\delta_1$, obtaining 
$\delta_1 \approx \deltabf \frac{\alpha'_{n}}{\alpha'_{n-1}}$.

To obtain a numerical approximation of the values $\alpha'_n$ 
we have computed points in the curve $(\alpha_n(\eps),\eps)\in 
S^-_n$ for small values of $\eps=2^{-k} h_n$ for different values 
of $k=1,2,..., M$ and $h_n$ a prescribed value depending on $n$ 
(which have been decreased when we increased $n$). 
Then, we have used this set of points to estimate the value 
$\alpha_n'\approx \frac{\alpha_n(2^{-k} h_n) -\alpha_n(0)}{2^{-k}h_n}$, 
and we did three steps of extrapolation to improve the accuracy of the 
results.  We 
have also considered consecutive approximations to estimate 
the accuracy of the results. These computations have been done with 
quadruple precision (using the library \cite{BHLR}). 

To compute each of the values $(\alpha_n(\eps),\eps)$ with good 
accuracy we have used the following procedure. We have computed
the $2^n$ periodic invariant curve of the system approximating 
it by its Fourier expansion. We use the same method described 
in Section 3.2.1 in \cite{JRT11p} (originally from \cite{CJ00}) 
to continue the zero Lyapunov exponent but tuned to 
continue the reducibility loss bifurcation curve. 
To continue this curve we use the characterization of the 
reducibility loss bifurcation given by Definition 2.3 in 
\cite{JRT11p}. 

\section{Description of the results} 

The values of $\alpha'_n$ actually depend on  $\omega$
(the rotation number of the system). Hence, from now on, we write 
$\alpha'_n=\alpha'_n(\omega)$. The estimated 
values of $\alpha'_n(\omega)$ 
for the family (\ref{FLM}) when $\omega= \frac{\sqrt{5}-1}{2}$ 
are shown in Table \ref{taula A}. We have 
included the ratios $\alpha'_n(\omega)/\alpha'_{n-1}(\omega)$ 
in the third column of the table. The value $\epsilon_a$ in the fourth 
column corresponds to the estimated accuracy (in absolute terms) for 
the value of $\alpha_n'$. In Tables \ref{taula 2} and \ref{taula 4} 
we included the same values for $\omega= 2 \frac{\sqrt{5}-1}{2}$ 
and $\omega= 4 \frac{\sqrt{5}-1}{2}$. 

We have computed the approximate value of $\beta'_i$ for the 
same values of $\omega$, but 
in all cases the value obtained has been equal (in the accuracy 
of the computations) to $-\alpha'_i$. Because of that 
these results will be omitted of the discussion.

We can observe in Tables \ref{taula A}, \ref{taula 2} and 
\ref{taula 4} that the ratios $\alpha'_i(\omega)/\alpha'_{i-1}(\omega)$ 
do not converge to a constant. Therefore (a priory) it seems that
there are no self-similarity of the bifurcation 
diagram. In fact, we will see in Section \ref{sec:Self renor bif diag} 
that one has to do a more subtle analysis to uncover the self-similarity 
properties of the family.

\subsection{Universality of the ratio sequence} 
\label{sec:Universality ratio sequence}

A remarkable fact is that the ratios $\alpha'_i(\omega)/\alpha'_{i-1}(\omega)$ 
on the third column of Table \ref{taula A} are approximately 
the same values on the third column of Table \ref{taula A}, but 
shifted on the index value $n$ by one position. 
Actually, the bigger is $n$ the closer are the values. 
The same phenomenon can be observed in Tables \ref{taula 2} 
and \ref{taula 4}.

Let us introduce some additional notation to follow with the analysis 
of this phenomenon. Consider $F$ a quasi-periodic forced map as follows,
\begin{equation}
\label{q.p. forced system}
\begin{array}{rccc}
F:& \T\times R &\rightarrow & \T \times R \\
  & \left( \begin{array}{c} \theta \\ x \end{array}\right)
  & \mapsto
  & \left( \begin{array}{c} \theta + \omega \\ f(\theta,x) \end{array}\right),
\end{array} 
\end{equation}
where $f$ is a $C^r$ map and $\omega$ an irrational number. This 
map $F$ can be identified with a pair $(\omega,f) \in \T\cap (\R\setminus\Q) 
 \times C^r(\T\times \R ,\R)$. Consider $\alpha_n'(\omega)$ 
the slope of the reducibility-loss bifurcation 
introduced before. Nothing  ensures yet their existence for a general 
map $F=(\omega,f)$, but whenever they exist they can be thought 
also depending, not only on $\omega$, but also on the function $f$ 
of the map $F$ considered \cite{fn1}.
In other words we consider $\alpha_n'(\omega) 
= \alpha_n'(\omega,f)$.

Given two sequences $\{r_i\}_{i\in \Z_+}$ and $\{s_i\}_{i\in \Z_+}$ of 
real numbers, we say that they are {\bf asymptotically 
equivalent} \cite{fn2}
if there exists $0<\rho <1$ and $K_0$ such that 
\[
| r_i - s_i| \leq K_0 \rho^i \quad \forall i\in \Z_+ .
\]
We denote this equivalence relation simply by $s_i \sim r_i$. 

Now we can resume the discussion on the phenomenon observed on Tables 
\ref{taula A}, \ref{taula 2} and \ref{taula 4}. Our numerical 
computations suggest that 
\begin{equation}
\label{equation equivalence of sequences} 
 \frac{\alpha_i'(\omega,f)}{\alpha_{i-1}'(\omega,f)} 
\sim \frac{\alpha_{i-1}'(2\omega,f)}{\alpha_{i-2}'(2\omega,f)}.
\end{equation}

Consider two different values $\omega_0$ and $\omega_1 $
 such that $2\omega_0 = 2 \omega_1 (\mod 1)$ . If equation 
(\ref{equation equivalence of sequences}) is true, then 
$\alpha_i'(\omega_0,f)/\alpha_{i-1}'(\omega_0,f)$ and 
$\alpha_i'(\omega_1,f)/\alpha_{i-1}'(\omega_1,f)$ should be 
asymptotically equivalent. This can be checked numerically. 
In Table \ref{taula B} we have recomputed the same values of 
Table \ref{taula A} but this time for $\omega = \frac{\sqrt{5}}{2}$ 
(and $f$ the function associated to the map (\ref{FLM})  as before). 
Again the results obtained show asymptotic equivalence.

\begin{table}[t!]
\begin{center}
\texttt{\begin{tabular}{|c|c|c|c|}
\hline 
\rule{0pt}{2.5ex} n & $\alpha_n'(\omega)$  &  $\alpha_n'(\omega)/\alpha_{n-1}'(\omega)$ &  $\epsilon_a$ \\ \hline 
0 &  -2.0000000000e+00 &  - - - &  3.0e-15 \\ 
1 &  -5.8329149229e+00 & 2.9164574614e+00 & 1.5e-14 \\ 
2 &  -8.4942599432e+00 & 1.4562633015e+00 & 8.4e-13 \\ 
3 &  -1.6351279467e+01 & 1.9249798777e+00 & 7.4e-15 \\ 
4 &  -1.1252460775e+01 & 6.8817004793e-01 & 3.0e-14 \\ 
5 &  -1.2243326651e+01 & 1.0880577054e+00 & 1.6e-13 \\ 
6 &  -1.8079693906e+01 & 1.4766978307e+00 & 1.6e-11 \\ 
7 &  -3.4735234067e+01 & 1.9212291009e+00 & 2.0e-12 \\ 
8 &  -2.9583312211e+01 & 8.5168023205e-01 & 2.1e-12 \\ 
9 &  -4.1569457725e+01 & 1.4051657715e+00 & 4.2e-10 \\ 
10 &  -7.8965495522e+01 & 1.8996036957e+00 & 9.1e-11 \\ 
11 &  -7.4500733455e+01 & 9.4345932945e-01 & 8.1e-10 \\ 
\hline
\end{tabular}}
\end{center}
\vspace{-4mm}
\caption{Approximate values of $\alpha_n'(\omega)$ for the
map (\ref{FLM}) for $ \omega = \frac{\sqrt{5} -1}{2}$.}
\label{taula A}

\vspace{2mm}

\begin{center}
\texttt{\begin{tabular}{|c|c|c|c|}
\hline 
 \rule{0pt}{2.5ex}  n & $\alpha_n'(\omega)$  &  $\alpha_n'(\omega)/\alpha_{n-1}'
(\omega)$ &  $\epsilon_a$  \\  \hline
0 &  -2.0000000000e+00 &  - - - &  4.1e-16 \\ 
1 &  -4.7793787548e+00 & 2.3896893774e+00 & 9.2e-14 \\ 
2 &  -9.9177338359e+00 & 2.0751094117e+00 & 8.0e-16 \\ 
3 &  -6.9333908531e+00 & 6.9909023249e-01 & 4.3e-15 \\ 
4 &  -7.5678156188e+00 & 1.0915028129e+00 & 2.4e-14 \\ 
5 &  -1.1183261803e+01 & 1.4777397292e+00 & 2.3e-12 \\ 
6 &  -2.1488744556e+01 & 1.9215095679e+00 & 3.0e-13 \\ 
7 &  -1.8302110429e+01 & 8.5170682641e-01 & 3.1e-13 \\ 
8 &  -2.5717669657e+01 & 1.4051750893e+00 & 6.1e-11 \\ 
9 &  -4.8853450105e+01 & 1.8996064090e+00 & 2.6e-11 \\ 
10 &  -4.6091257360e+01 & 9.4345961772e-01 & 1.2e-10 \\ 
11 &  -7.1498516059e+01 & 1.5512381339e+00 & 4.1e-09 \\ 
\hline
\end{tabular}}
\end{center}
\vspace{-4mm}
\caption{Approximate values of $\alpha_n'(\omega)$ for the
map (\ref{FLM}) for $ \omega = 2 \frac{\sqrt{5} -1}{2}$.}
\label{taula 2}

\vspace{2mm}

\begin{center}
\texttt{\begin{tabular}{|c|c|c|c|}
\hline 
 \rule{0pt}{2.5ex}  n & $\alpha_n'(\omega)$  &  $\alpha_n'(\omega)/\alpha_{n-1}'
(\omega)$ &  $\epsilon_a$  \\  \hline
0 &  -2.0000000000e+00 &  - - - &  5.3e-15 \\ 
1 &  -6.1135714539e+00 & 3.0567857269e+00 & 4.0e-17 \\ 
2 &  -4.6432689399e+00 & 7.5950186809e-01 & 8.8e-16 \\ 
3 &  -5.1662366620e+00 & 1.1126292121e+00 & 5.4e-15 \\ 
4 &  -7.6637702641e+00 & 1.4834338350e+00 & 5.2e-13 \\ 
5 &  -1.4738755184e+01 & 1.9231728870e+00 & 6.7e-14 \\ 
6 &  -1.2555429245e+01 & 8.5186497017e-01 & 6.9e-14 \\ 
7 &  -1.7643286059e+01 & 1.4052316105e+00 & 1.3e-11 \\ 
8 &  -3.3515585777e+01 & 1.8996226477e+00 & 5.8e-12 \\ 
9 &  -3.1620659727e+01 & 9.4346134772e-01 & 2.6e-11 \\ 
10 &  -4.9051192417e+01 & 1.5512387420e+00 & 9.0e-10 \\ 
11 &  -9.2911119039e+01 & 1.8941663691e+00 & 5.6e-12 \\ 
\hline
\end{tabular}}
\end{center}
\vspace{-4mm}
\caption{Approximate values of $\alpha_n'(\omega)$ for the
map (\ref{FLM}) for $ \omega = 4 \frac{\sqrt{5} -1}{2}$.}
\label{taula 4}
\end{table}

\begin{table}[t!]
\begin{center}
\texttt{\begin{tabular}{|c|c|c|c|}
\hline 
 \rule{0pt}{2.5ex}  n & $\alpha_n'(\omega)$  &  $\alpha_n'(\omega)/\alpha_{n-1}'(\omega)$ &  $\epsilon_a$ \\ \hline 
0 &  -2.0000000000e+00 &  - - - &  5.3e-17 \\ 
1 &  -3.7459592187e+00 & 1.8729796094e+00 & 3.6e-15 \\ 
2 &  -6.1798012892e+00 & 1.6497246575e+00 & 2.4e-13 \\ 
3 &  -1.2170152952e+01 & 1.9693437350e+00 & 2.1e-15 \\ 
4 &  -8.4095313813e+00 & 6.9099635922e-01 & 9.5e-15 \\ 
5 &  -9.1576001570e+00 & 1.0889548706e+00 & 5.2e-14 \\ 
6 &  -1.3525371011e+01 & 1.4769558377e+00 & 5.0e-12 \\ 
7 &  -2.5986306241e+01 & 1.9213008072e+00 & 6.5e-13 \\ 
8 &  -2.2132200091e+01 & 8.5168703416e-01 & 6.6e-13 \\ 
9 &  -3.1099463199e+01 & 1.4051681745e+00 & 1.3e-10 \\ 
10 &  -5.9076676864e+01 & 1.8996043914e+00 & 5.6e-11 \\ 
11 &  -5.5736446311e+01 & 9.4345940345e-01 & 2.5e-10 \\ 
\hline
\end{tabular}}
\end{center}
\vspace{-4mm}
\caption{Approximate values of $\alpha_n'(\omega)$ for the
map (\ref{FLM}) for $ \omega =  \frac{\sqrt{5}}{2}$.}

\label{taula B}

\vspace{2mm}

\begin{center}
\texttt{\begin{tabular}{|c|c|c|c|}
\hline 
 \rule{0pt}{2.5ex}  n & $\alpha_n'(\omega)$  &  $\alpha_n'(\omega)/\alpha_{n-1}'(\omega)$ &  $\epsilon_a$ \\ \hline
0 &  -4.0000000000e+00 &  - - - &  2.7e-14 \\ 
1 &  -8.1607837043e+00 & 2.0401959261e+00 & 5.1e-14 \\ 
2 &  -1.1166652707e+01 & 1.3683309241e+00 & 2.4e-12 \\ 
3 &  -2.1221554117e+01 & 1.9004400578e+00 & 2.1e-14 \\ 
4 &  -1.4564213015e+01 & 6.8629342294e-01 & 8.6e-14 \\ 
5 &  -1.5837452605e+01 & 1.0874224778e+00 & 4.6e-13 \\ 
6 &  -2.3384207858e+01 & 1.4765132021e+00 & 4.5e-11 \\ 
7 &  -4.4925217655e+01 & 1.9211776567e+00 & 5.8e-12 \\ 
8 &  -3.8261700375e+01 & 8.5167534788e-01 & 5.9e-12 \\ 
9 &  -5.3763965691e+01 & 1.4051640456e+00 & 1.1e-09 \\ 
10 &  -1.0213020106e+02 & 1.8996031960e+00 & 1.1e-10 \\ 
11 &  -9.6355685578e+01 & 9.4345927630e-01 & 2.3e-09 \\ 
\hline
\end{tabular}}
\end{center}
\vspace{-4mm}
\caption{Approximate values of $\alpha_n'$ for the
map (\ref{FLM B}) for $ \omega =  \frac{\sqrt{5} -1}{2}$.}
\label{taula ABis}

\vspace{2mm}

\begin{center}
\texttt{\begin{tabular}{|c|c|c|c|}
\hline 
 \rule{0pt}{2.5ex} n & $\alpha_n'(\omega)$  &  $\alpha_n'(\omega)/\alpha_{n-1}'(\omega)$ &  $\epsilon_a$ \\ \hline 
0 &  -4.0000000000e+00 &  - - - &  3.8e-15 \\ 
1 &  -6.2417012728e+00 & 1.5604253182e+00 & 2.4e-13 \\ 
2 &  -1.2036825830e+01 & 1.9284527253e+00 & 2.1e-15 \\ 
3 &  -8.2891818940e+00 & 6.8865180997e-01 & 9.0e-15 \\ 
4 &  -9.0187641307e+00 & 1.0880161934e+00 & 4.9e-14 \\ 
5 &  -1.3318271659e+01 & 1.4767291245e+00 & 4.7e-12 \\ 
6 &  -2.5587438370e+01 & 1.9212281462e+00 & 6.1e-13 \\ 
7 &  -2.1792312360e+01 & 8.5168011135e-01 & 6.2e-13 \\ 
8 &  -3.0621808757e+01 & 1.4051656497e+00 & 1.2e-10 \\ 
9 &  -5.8169300479e+01 & 1.8996036760e+00 & 5.2e-11 \\ 
10 &  -5.4880369083e+01 & 9.4345932701e-01 & 2.3e-10 \\ 
11 &  -8.5132515708e+01 & 1.5512380316e+00 & 8.2e-09 \\ 
\hline
\end{tabular}}
\end{center}
\vspace{-4mm}
\caption{Approximate values of $\alpha_n'$ for the
map (\ref{FLM B}) for $ \omega = 2 \frac{\sqrt{5} -1}{2}$.}
\label{taula 2Bis}
\end{table}

Given $F$ a skew product map like (\ref{q.p. forced system}), let $F^2$ denote 
the map composed with itself. Concretely, the map $F^2$ is given as 
$F^2(\theta, x) = (\theta + 2 \omega, f(\theta + \omega, f(\theta,x)))$.
To simplify the notation, let us denote by $f^2$ the map defined as 
$f^2(\theta,x):= f(\theta + \omega, f(\theta,x))$. Then we have that 
the map $F^2$ corresponds to the pair $(2\omega, f^2)$. Moreover we have 
that the $2^k$ periodic curves
are $2^{k-1}$ periodic curves of the map $F^2$, therefore
$\alpha_i'(\omega,f) = \alpha_{i-1}'(2\omega, f^2)$. Hence the 
relation given by (\ref{equation equivalence of sequences}) can be rewritten 
as 
\begin{equation}
\label{equation equivalence of sequences 2}
\frac{\alpha_i'(2\omega,f^2)}{\alpha_{i-1}'(2\omega,f^2)} 
\sim \frac{\alpha_{i}'(2\omega,f)}{\alpha_{i-1}'(2\omega,f)} .
\end{equation}

Now we have that the relation given by (\ref{equation equivalence of sequences}) can be 
explained as a consequence of a much more general phenomenon. 
We believe that the sequence 
$\frac{\alpha_i'(2\omega,f)}{\alpha_{i-1}'(2\omega,f)}$ is (asymptotically) {\bf universal}, 
in the sense that it does not depend on the map $f$. This would 
imply (\ref{equation equivalence of sequences 2}) and consequently 
(\ref{equation equivalence of sequences}).

To check the universality of the sequence, we consider the following map,
\begin{equation}
\label{FLM B}
\left.
\begin{array}{rcl}
\bar{\theta} & = &  \theta + \omega  ,\\
\bar{x} & = & \alpha x(1-x) +  \eps \cos( 2\pi\theta) .
\end{array}
\right\}
\end{equation}
This map is like (\ref{FLM}) but with an additive forcing
instead of a multiplicative one. In the literature, sometimes (\ref{FLM}) is 
referred as the Driven Logistic Map and (\ref{FLM B}) is referred as the 
Forced Logistic Map. We do not do this distinction in this paper and 
we consider both as two different versions of the Forced Logistic Map. 
 
Note that both maps are in the class of quasi-periodically 
forced one dimensional unimodal maps, with a quasi-periodic forcing
of the type $h(x)\cos(2\pi\theta)$, where $h$ is a function of one variable.
This is certainly a very restrictive class of maps. In Section 
\ref{subsection counter example of universality.}
we explore what happens when the quasi-periodic forcing is not of 
this form. 

We have computed the slopes $\alpha_n'(\omega)$, the associated 
ratios  $\alpha_n'(\omega)/ \alpha_{n-1}'(\omega)$
and the estimation of the accuracy for the 
family (\ref{FLM B}) as we did before for the family (\ref{FLM}). 
The results are shown in Table \ref{taula ABis} for 
$\omega=\frac{\sqrt{5} -1}{2}$ and Table \ref{taula 2Bis} 
for $\omega=2\frac{\sqrt{5} -1}{2}$.

Now we can compare the sequences  $\alpha_i'(\omega_0)/\alpha_{i-1}'(\omega_0)$ 
in Table \ref{taula A} (respectively \ref{taula 2}) with the ones of 
Table \ref{taula ABis} (respectively \ref{taula 2Bis}).
Again we can observe that both sequences have
an equivalent asymptotic behavior for equal values of $\omega$.
This agrees with the conjectured universal behavior. 

\subsection{Self-similarity of the bifurcation diagram}
\label{sec:Self renor bif diag}

\begin{table}[t]
\begin{center}
\texttt{\begin{tabular}{|c|c|c||c|c|}
\hline 
 n & $ \delta_{1,n} (\omega_0)$  & $\delta_{1,n}  - \delta_{1,n-1} $  &  $ \delta_{1,n} (\omega_0) $ & $\delta_{1,n} - \delta_{1,n-1} $ \\ \hline
1 &  1.36175279e+01 &    - - -  & 1.11579415e+01 &    - - -  \\ 
2 &  8.29844510e+00 &   -5.3e+00 & 7.57460662e+00 &   -3.6e+00 \\ 
3 &  7.69807112e+00 &   -6.0e-01 & 6.97211386e+00 &   -6.0e-01 \\ 
4 &  7.57782290e+00 &   -1.2e-01 & 6.83972864e+00 &   -1.3e-01 \\ 
5 &  7.55390503e+00 &   -2.4e-02 & 6.81347460e+00 &   -2.6e-02 \\ 
6 &  7.54857906e+00 &   -5.3e-03 & 6.80758174e+00 &   -5.9e-03 \\ 
7 &  7.54747726e+00 &   -1.1e-03 & 6.80631795e+00 &   -1.3e-03 \\ 
8 &  7.54724159e+00 &   -2.4e-04 & 6.80604419e+00 &   -2.7e-04 \\ 
9 &  7.54719154e+00 &   -5.0e-05 & 6.80598601e+00 &   -5.8e-05 \\ 
10 &  7.54718076e+00 &   -1.1e-05 & 6.80597353e+00 &   -1.2e-05 \\ 
11 &  7.54717846e+00 &   -2.3e-06 & 6.80597086e+00 &   -2.7e-06 \\ 
\hline
\end{tabular}}
\end{center}
\vspace{-4mm}
\caption{ 
Estimations of the value $\delta_{1,n}= \delta_{1,n}(\omega) 
=\deltabf \frac{\alpha_n'(\omega_0)}{\alpha_{n-1}'(2 \omega_0)}$
for $\omega_0 = \frac{\sqrt{5}-1}{2}$ on the left
and $\omega_1= 2 \frac{\sqrt{5}-1}{2}$ on the right. }
\label{taula relacions}
\end{table}

\begin{table}[t]
\begin{center}
\texttt{\begin{tabular}{|c|c|c|}
\hline 
 n & $ \delta_{1,n} (\omega_0)$ & $\delta_{1,n} - \delta_{1,n-1} $ \\ \hline
1 &  9.52608610e+00 &    - - -  \\ 
2 &  8.35338804e+00 &   -1.2e+00 \\ 
3 &  8.23204689e+00 &   -1.2e-01 \\ 
4 &  8.20385506e+00 &   -2.8e-02 \\ 
5 &  8.19937833e+00 &   -4.5e-03 \\ 
6 &  8.19817945e+00 &   -1.2e-03 \\ 
7 &  8.19796400e+00 &   -2.2e-04 \\ 
8 &  8.19791815e+00 &   -4.6e-05 \\ 
9 &  8.19790879e+00 &   -9.4e-06 \\ 
10 &  8.19790672e+00 &   -2.1e-06 \\ 
11 &  8.19790628e+00 &   -4.4e-07 \\ 
\hline
\end{tabular}}
\end{center}
\vspace{-4mm}
\caption{ 
Estimations of the value
$\delta_{1,n}= \delta_{1,n}(\omega) =\deltabf \frac{\alpha_n'(\omega_0)}{\alpha_{n-1}'(2\omega_0)}$
for $\omega_0 = \frac{\sqrt{5}-1}{2}$ for the family (\ref{FLM B}) 
 with $\omega_0 = \frac{\sqrt{5}-1}{2}$ . }
\label{taula relacions Bis}
\end{table}

In order to find self-similarity properties of the parameter space 
we need to refine our analysis. Given a one dimensional 
map in the interval  $g:I\rightarrow I$,  its (doubling) renormalization is
defined as $\RR(f) = A^{-1}\circ g \circ g \circ A$ with $A$ an
affine transformation.
Given a q.p. forced map like (\ref{q.p. forced system}),
one might try to define a renormalization for these kind of maps. 
The most simple choice would be to define  its
(doubling) renormalization  also 
as $\TT(F) = A^{-1}\circ F \circ F \circ A$, 
with $A$ a suitable affine map. If the map $F$ has rotation
number $\omega$, its (doubling) renormalization will
have rotation number $2\omega$.
This is due to the composition of $F$ with itself when we
define its renormalization.
This argument  become more clear in the rigorous definition of
the renormalization operator for quasi-periodically maps done in 
\cite{JRT11a}. 

\begin{figure}[t]
\begin{center}
\includegraphics[width=15cm]{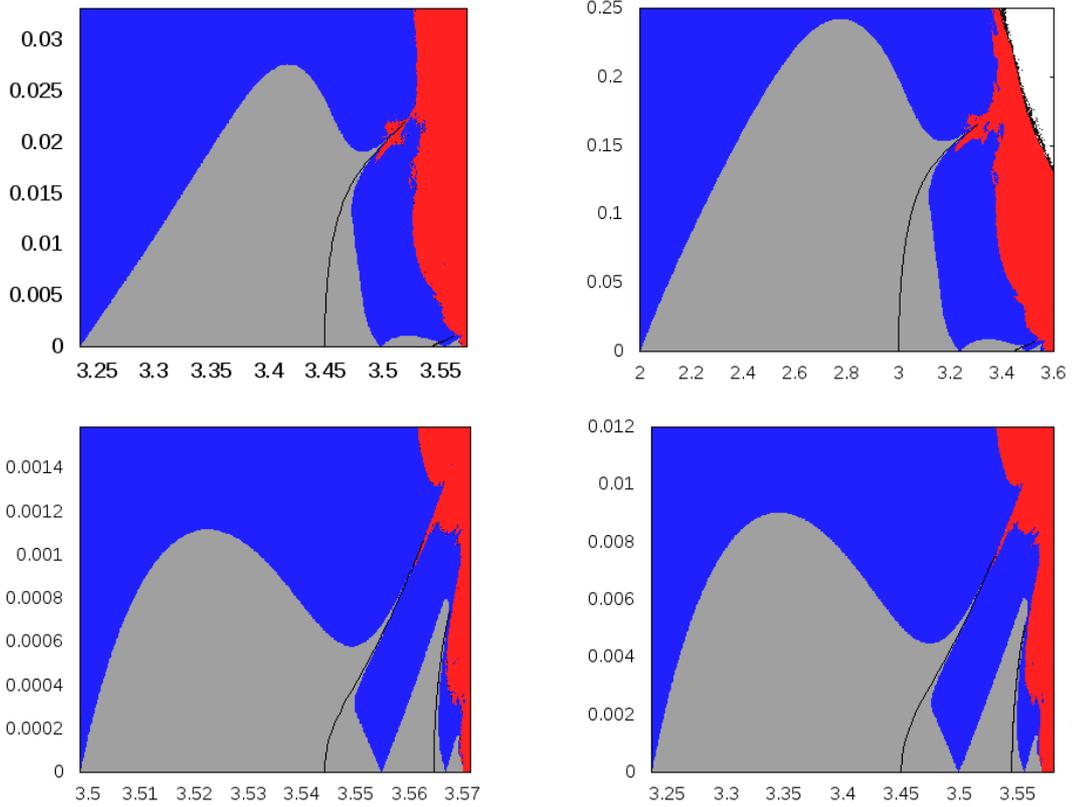}
\caption{  
Diagram of the parameter space of the map (\ref{FLM})
for  $\omega= \frac{\sqrt{5}-1}{2}$ on the left hand side 
and for $\omega= \sqrt{5}-1$ on the right. See Table 
\ref{table coding diagrams} for the
coding of each color and the text for more details.}
\label{FLM parameter space 2}
\end{center}
\end{figure}

Consider the reducibility loss curves $S^\pm_n$ as 
in the parameter space of the map (\ref{FLM})
(see Section \ref{Sec:Descrip Comp}). 
We have that these curves depend on the rotation number 
(i.e. $S^\pm_{n}=S^\pm_{n}(\omega)$). Now we can look 
for self-similarity of the parameter space as 
in Section \ref{Sec:Descrip Comp} but  taking into account 
the doubling of the rotation number. In other words we can 
look for affine relationship between the bifurcation diagram 
around $s_n$ for $\omega_0$ and the bifurcation diagram 
around $s_{n-1}$ for $2 \omega_0$. Consider an affine map like 
(\ref{equation afine transformation parameters}), but such that 
it maps the curves $S^-_{i}(\omega)$ to $S^-_{i-1}(2\omega)$
(and respectively $S^+_{i}(\omega)$ to $S^+_{i-1}(2\omega)$).
Redoing the same  computation of Section \ref{Sec:Descrip Comp} 
we have that if there exists an affine relation between the parameters spaces it
should be given by $\delta_1 \approx \deltabf  
\frac{\alpha_i(\omega)}{\alpha_{i-1}(2\omega)}$.

In Table \ref{taula relacions} we have estimations of
these values for the family (\ref{FLM}) for
$\omega= \frac{\sqrt{5} -1}{2}$ (left) and $\omega= \sqrt{5} -1$ (right).
In Table \ref{taula relacions Bis} we have the same estimation 
for  the family (\ref{FLM B}) and $\omega= \frac{\sqrt{5} -1}{2}$. 
In all cases the sequence converges to a concrete value. This convergence 
means that there exists an affine relation between the reducibility 
loss bifurcation curves around the superstable  periodic
orbits of the uncoupled map. Note that the limit constants  
obtained in each of these three cases are different one to each 
other. This indicates that the renormalization factor depends 
both on the value of $\omega$ taken and the family of maps considered. 
Therefore the renormalization factor 
$\delta_1$ is not universal.

These evidences of self-similarity are only valid for 
infinitely small values of the coupling parameter $\eps$. 
If the self-similarity between families extends to larger values of $\eps$, 
the bifurcation diagram of the map (\ref{FLM}) in a box $I_1\times I_2$ 
of the parameters space for a prescribed $\omega=\omega_0$ should 
be approximately the same as the diagram in the 
box $L(I_1\times I_2)$ for $\omega=2\omega_0$, where $L$ is the 
affine map (\ref{equation afine transformation parameters}).

In Figure  \ref{FLM parameter space 2} we have 
a bifurcation diagram of the map (\ref{FLM}) analog 
to the one displayed in Figure \ref{FLM parameter space}. 
The boxes 
have been selected such that the one in the
left is the image of the one in the right through the affine map
$L$ given by (\ref{equation afine transformation parameters}).
The value of $\delta_0$ has been taken equal to
the Feigenbaum constant ($\approx \texttt{4,66920}$) and  
$\delta_1 \approx \texttt{7,54718}$ the experimental value 
obtained in Table \ref{taula relacions}. The results indicate
that self-similarity properties extend to the whole reducibility 
region around each period doubling bifurcation. 

\subsection{Non-universality of the ratio sequence under 
more general conditions}
\label{subsection counter example of universality.}

\begin{table}[t!]
\begin{center}
\texttt{\begin{tabular}{|c|c|c|c|}
\hline 
 n & $\alpha_n'(\omega)$  &  $\alpha_n'(\omega)/\alpha_{n-1}'(\omega)$ &  $\epsilon_a$ \\ \hline 
0 &  -4.4000000000e+00 &  - - - &  2.1e-14 \\ 
1 &  -8.5708073961e+00 & 1.9479107718e+00 & 3.2e-14 \\ 
2 &  -1.2367363641e+01 & 1.4429636637e+00 & 2.5e-12 \\ 
3 &  -2.2002414361e+01 & 1.7790707058e+00 & 2.0e-14 \\ 
4 &  -1.5466051366e+01 & 7.0292519321e-01 & 1.2e-13 \\ 
5 &  -1.7124001858e+01 & 1.1071993396e+00 & 7.5e-13 \\ 
6 &  -2.5233583736e+01 & 1.4735798293e+00 & 5.0e-11 \\ 
7 &  -4.5526415150e+01 & 1.8041993411e+00 & 4.4e-12 \\ 
8 &  -4.0793050977e+01 & 8.9603037802e-01 & 3.7e-12 \\ 
9 &  -5.9579098646e+01 & 1.4605207804e+00 & 1.1e-09 \\ 
10 &  -1.0695246126e+02 & 1.7951339260e+00 & 9.0e-11 \\ 
11 &  -1.0464907069e+02 & 9.7846341692e-01 & 3.1e-09 \\ 
\hline
\end{tabular}}
\end{center}
\vspace{-4mm}
\caption{ 
Approximate values of $\alpha_n'(\omega)$ of the 
family (\ref{FLM C}) for $ \omega = \frac{\sqrt{5} -1}{2}$ and $E=10^{-1}$.}

\label{taula contraexemple 1.1} 

\vspace{1mm}

\begin{center}
\texttt{\begin{tabular}{|c|c|c|c|}
\hline 
 n & $\alpha_n'(\omega)$  &  $\alpha_n'(\omega)/\alpha_{n-1}'(\omega)$ &  $\epsilon_a$ \\ \hline
0 &  -4.4000000000e+00 &  - - - &  4.2e-15 \\ 
1 &  -6.8849393233e+00 & 1.5647589371e+00 & 2.6e-13 \\ 
2 &  -1.2339134533e+01 & 1.7921921972e+00 & 1.9e-15 \\ 
3 &  -8.8252833730e+00 & 7.1522709709e-01 & 1.0e-14 \\ 
4 &  -9.9507316299e+00 & 1.1275254526e+00 & 2.6e-14 \\ 
5 &  -1.5240261120e+01 & 1.5315719172e+00 & 4.9e-12 \\ 
6 &  -2.6987715314e+01 & 1.7708171207e+00 & 6.9e-13 \\ 
7 &  -2.4101506278e+01 & 8.9305471016e-01 & 1.0e-12 \\ 
8 &  -3.4274713216e+01 & 1.4220983876e+00 & 1.4e-10 \\ 
9 &  -5.9739968275e+01 & 1.7429750002e+00 & 4.0e-11 \\ 
10 &  -6.0772668123e+01 & 1.0172865818e+00 & 3.5e-10 \\ 
11 &  -9.7063373626e+01 & 1.5971550472e+00 & 8.3e-09 \\ 
\hline
\end{tabular}}
\end{center}
\vspace{-3mm}
\caption{ 
Approximate values of $\alpha_n'(\omega)$ of the 
family (\ref{FLM C}) for $ \omega = 2\frac{\sqrt{5} -1}{2}$ and $E=10^{-1}$.}
\label{taula contraexemple 1.2}

\vspace{1mm}

\begin{center}
\texttt{\begin{tabular}{|c|c|c|c|}
\hline 
 n & $\alpha_n'(\omega)$  &  $\alpha_n'(\omega)/\alpha_{n-1}'(\omega)$ &  $\epsilon_a$\\ \hline 
0 &  -4.0400000000e+00 &  - - - &  2.6e-14 \\ 
1 &  -8.1970065912e+00 & 2.0289620275e+00 & 4.9e-14 \\ 
2 &  -1.1286613379e+01 & 1.3769189098e+00 & 2.4e-12 \\ 
3 &  -2.1298993542e+01 & 1.8871022535e+00 & 2.1e-14 \\ 
4 &  -1.4654396022e+01 & 6.8803232383e-01 & 8.9e-14 \\ 
5 &  -1.5964802815e+01 & 1.0894207302e+00 & 5.0e-13 \\ 
6 &  -2.3545880397e+01 & 1.4748619617e+00 & 4.5e-11 \\ 
7 &  -4.4967877741e+01 & 1.9097981041e+00 & 5.6e-12 \\ 
8 &  -3.8501241787e+01 & 8.5619432630e-01 & 4.8e-12 \\ 
9 &  -5.4345411904e+01 & 1.4115236128e+00 & 1.1e-09 \\ 
10 &  -1.0260055779e+02 & 1.8879341272e+00 & 1.0e-10 \\ 
11 &  -9.7178352178e+01 & 9.4715227943e-01 & 2.3e-09 \\ 
\hline
\end{tabular}}
\end{center}
\vspace{-3mm}
\caption{ 
Approximate values of $\alpha_n'(\omega)$ of the 
family (\ref{FLM C}) for $ \omega =  \frac{\sqrt{5} -1}{2}$ and $E=10^{-2}$.}
\label{taula contraexemple 01.1}
\end{table}

\begin{table}[t!]
\begin{center}
\texttt{\begin{tabular}{|c|c|c|c|}
\hline 
 n & $\alpha_n'(\omega)$  &  $\alpha_n'(\omega)/\alpha_{n-1}'(\omega)$ &  $\epsilon_a$ \\ \hline  
0 &  -4.0400000000e+00 &  - - - &  3.8e-15 \\ 
1 &  -6.2967793362e+00 & 1.5586087466e+00 & 2.4e-13 \\ 
2 &  -1.2062745707e+01 & 1.9157008786e+00 & 2.1e-15 \\ 
3 &  -8.3390918858e+00 & 6.9130959805e-01 & 9.0e-15 \\ 
4 &  -9.1092083960e+00 & 1.0923501648e+00 & 4.6e-14 \\ 
5 &  -1.3509885744e+01 & 1.4831020608e+00 & 4.7e-12 \\ 
6 &  -2.5722289496e+01 & 1.9039605503e+00 & 6.1e-13 \\ 
7 &  -2.2023032848e+01 & 8.5618478288e-01 & 6.8e-13 \\ 
8 &  -3.0953588852e+01 & 1.4055098163e+00 & 1.2e-10 \\ 
9 &  -5.8281365651e+01 & 1.8828629510e+00 & 5.1e-11 \\ 
10 &  -5.5448581611e+01 & 9.5139468666e-01 & 2.5e-10 \\ 
11 &  -8.6313143445e+01 & 1.5566339289e+00 & 8.2e-09 \\ 
\hline
\end{tabular}}
\end{center}
\vspace{-3mm}
\caption{
Approximate values of $\alpha_n'(\omega)$ of the family 
(\ref{FLM C}) for $ \omega =2 \frac{\sqrt{5} -1}{2}$ and $E=10^{-2}$.}
\label{taula contraexemple 01.2} 

\vspace{1mm}

\begin{center}
\texttt{\begin{tabular}{|c|c|c|c|}
\hline 
 n & $\alpha_n'(\omega)$  &  $\alpha_n'(\omega)/\alpha_{n-1}'(\omega)$ &  $\epsilon_a$ \\ \hline 
0 &  -4.0040000000e+00 &  - - - &  2.7e-14 \\ 
1 &  -8.1643491058e+00 & 2.0390482282e+00 & 5.0e-14 \\ 
2 &  -1.1178647202e+01 & 1.3692024995e+00 & 2.4e-12 \\ 
3 &  -2.1229290655e+01 & 1.8990930004e+00 & 2.1e-14 \\ 
4 &  -1.4573231305e+01 & 6.8646812284e-01 & 8.6e-14 \\ 
5 &  -1.5850170411e+01 & 1.0876222355e+00 & 4.7e-13 \\ 
6 &  -2.3400073191e+01 & 1.4763294391e+00 & 4.5e-11 \\ 
7 &  -4.4929300750e+01 & 1.9200495820e+00 & 5.7e-12 \\ 
8 &  -3.8285483755e+01 & 8.5212730036e-01 & 5.8e-12 \\ 
9 &  -5.3822109356e+01 & 1.4058098286e+00 & 1.1e-09 \\ 
10 &  -1.0217709608e+02 & 1.8984223641e+00 & 1.0e-10 \\ 
11 &  -9.6437862923e+01 & 9.4383053171e-01 & 2.2e-09 \\ 
\hline
\end{tabular}}
\end{center}
\vspace{-3mm}
\caption{
Approximate values of $\alpha_n'(\omega)$ of the family 
(\ref{FLM C}) for $ \omega = \frac{\sqrt{5} -1}{2}$ and $E=10^{-3}$.}
\label{taula contraexemple 001.1}

\vspace{1mm}

\begin{center}
\texttt{\begin{tabular}{|c|c|c|c|}
\hline 
 n & $\alpha_n'(\omega)$  &  $\alpha_n'(\omega)/\alpha_{n-1}'(\omega)$ &  $\epsilon_a$ \\ \hline
0 &  -4.0040000000e+00 &  - - - &  3.7e-15 \\ 
1 &  -6.2470810635e+00 & 1.5602100558e+00 & 2.4e-13 \\ 
2 &  -1.2039370639e+01 & 1.9271993618e+00 & 2.1e-15 \\ 
3 &  -8.2941266770e+00 & 6.8891696467e-01 & 9.0e-15 \\ 
4 &  -9.0277699390e+00 & 1.0884533466e+00 & 4.8e-14 \\ 
5 &  -1.3337423899e+01 & 1.4773774686e+00 & 4.7e-12 \\ 
6 &  -2.5600860673e+01 & 1.9194756698e+00 & 6.1e-13 \\ 
7 &  -2.1815381589e+01 & 8.5213469453e-01 & 6.3e-13 \\ 
8 &  -3.0654499563e+01 & 1.4051782426e+00 & 1.2e-10 \\ 
9 &  -5.8180013759e+01 & 1.8979273708e+00 & 5.2e-11 \\ 
10 &  -5.4936892373e+01 & 9.4425712239e-01 & 2.4e-10 \\ 
11 &  -8.5250385421e+01 & 1.5517875464e+00 & 8.2e-09 \\ 
\hline
\end{tabular}}
\end{center}
\vspace{-3mm}
\caption{ 
Approximate values of $\alpha_n'(\omega)$ of the family
(\ref{FLM C}) for $ \omega = 2\frac{\sqrt{5} -1}{2}$ and $E=10^{-3}$.}
\label{taula contraexemple 001.2}
\end{table}

\begin{table}[t!]
\begin{center}
\hspace{-4.5mm}
\mbox{
\texttt{\begin{tabular}{|c|c|c|c|}
\hline 
 n & $ \delta_{1,n} (\omega_0)$ & $\delta_{1,n} - \delta_{1,n-1} $ \\ \hline
1 &  9.095188 &    - - -  \\ 
2 &  8.387251 &   -7.1e-01 \\ 
3 &  8.325844 &   -6.1e-02 \\ 
4 &  8.182639 &   -1.4e-01 \\ 
5 &  8.035129 &   -1.5e-01 \\ 
6 &  7.730884 &   -3.0e-01 \\ 
7 &  7.876621 &   1.5e-01 \\ 
8 &  7.902866 &   2.6e-02 \\ 
9 &  8.116387 &   2.1e-01 \\ 
10 &  8.359271 &   2.4e-01 \\ 
11 &  8.040253 &   -3.2e-01 \\ 
\hline
\end{tabular}}
\texttt{\begin{tabular}{|c|c|c|c|}
\hline 
 n & $ \delta_{1,n} (\omega_0)$ & $\delta_{1,n} - \delta_{1,n-1} $ \\ \hline
1 &  9.473633 &    - - -  \\ 
2 &  8.369274 &   -1.1e+00 \\ 
3 &  8.244333 &   -1.2e-01 \\ 
4 &  8.205249 &   -3.9e-02 \\ 
5 &  8.183245 &   -2.2e-02 \\ 
6 &  8.137779 &   -4.5e-02 \\ 
7 &  8.162729 &   2.5e-02 \\ 
8 &  8.162820 &   9.1e-05 \\ 
9 &  8.197747 &   3.5e-02 \\ 
10 &  8.219826 &   2.2e-02 \\ 
11 &  8.183173 &   -3.7e-02 \\ 
\hline
\end{tabular}}
\texttt{\begin{tabular}{|c|c|c|c|}
\hline 
 n & $ \delta_{1,n} (\omega_0)$ & $\delta_{1,n} - \delta_{1,n-1} $ \\ \hline
1 &  9.520727 &    - - -  \\ 
2 &  8.355159 &   -1.2e+00 \\ 
3 &  8.233307 &   -1.2e-01 \\ 
4 &  8.204041 &   -2.9e-02 \\ 
5 &  8.197777 &   -6.3e-03 \\ 
6 &  8.191961 &   -5.8e-03 \\ 
7 &  8.194411 &   2.4e-03 \\ 
8 &  8.194339 &   -7.1e-05 \\ 
9 &  8.198023 &   3.7e-03 \\ 
10 &  8.200161 &   2.1e-03 \\ 
11 &  8.196456 &   -3.7e-03 \\ 
\hline
\end{tabular}}}
\end{center}
\vspace{-4mm}
\caption{ 
Estimations of the value  $\delta_{1,n}= \delta_{1,n}(\omega) 
=\deltabf \frac{\alpha_n'(\omega_0)}{\alpha_{n-1}'(2 \omega_0)}$
for the map (\ref{FLM C}) (where $\deltabf$ is the Feigenbaum constant). 
The different boxes correspond (from left to right) to $E=10^{-1},10^{-2}$
and $10^{-3}$. In all the cases 
we have taken $\omega=  \frac{\sqrt{5}-1}{2}$.}
\label{taula relacions afins contraexemple}
\end{table}

A natural question to ask after the numerical evidences of 
universality and renormalization reported in the previous 
sections is ``how general are these phenomena?''. In this 
section we present an example which demonstrates that the universality 
and the self similarity properties depend on the Fourier expansion 
of the quasi-periodic forcing. 
We must say that we designed this example after developing some
of the theory presented in \cite{JRT11a,JRT11b,JRT11c} and 
the cited theory provides a theoretical explanation to both 
behaviors.

Consider the following family, 
\begin{equation}
\label{FLM C}
\left.
\begin{array}{rcl}
\bar{\theta} & = & \theta + \omega  ,\\
\bar{x} & = & \alpha x(1-x) + \eps (\cos(2\pi \theta) + E\cos(4 \pi \theta)) .
\end{array}
\right\}
\end{equation}
We consider $E$ a fixed value and 
$\alpha$ and $\eps$ as true parameters, obtaining a two parametric family. 
Remark that for $E=0$ we recover the family  (\ref{FLM B}) introduced in 
Section \ref{sec:Universality ratio sequence}. 

We have done the computation of the values  $\alpha_n'(\omega)$  
for the family (\ref{FLM C}), for $E=10^{-1},10^{-2},10^{-3}$ and 
$\omega=\frac{\sqrt{5} -1}{2}, 2 \frac{\sqrt{5} -1}{2}$. 
The results are shown in Tables \ref{taula contraexemple 1.1} to
\ref{taula contraexemple 001.2}. 
To compute the values in the tables we have used the same procedure 
that we have used for the families (\ref{FLM}) and (\ref{FLM B}) before. 
In these tables we have also included
 the estimated values of the ratios $\alpha_n'(\omega)/\alpha_{n-1}'(\omega)$ 
and the estimated accuracies. In Table \ref{taula relacions afins contraexemple} 
we include the ratios $\delta_{1,n} (\omega) = 
\alpha_n'(\omega)/\alpha_{n-1}'(2\omega)$ and the differences 
 $\delta_{1,n} (\omega) -  \delta_{1,n-1} (\omega)$. 

In the third column of Tables \ref{taula contraexemple 1.1} 
to \ref{taula contraexemple 001.2} we can observe
that the sequence of ratios $\alpha_n'(\omega)/\alpha_{n-1}'(\omega)$ 
ceases to be universal. In other words, the sequence is not 
asymptotically equivalent to the sequences obtained for 
the maps (\ref{FLM}) and (\ref{FLM B}) (displayed in 
Tables \ref{taula A}  to \ref{taula 2Bis}). 
We can also observe in Table 
\ref{taula relacions afins contraexemple} that the different 
sequences $\deltabf \alpha_n'(\omega)/\alpha_{n-1}'(2\omega)$ 
cease to converge. Recall that the limit of this sequence gives 
us the scale factor between the bifurcations diagram of the map 
and itself for a doubled period. In other words, the self-similarity 
properties of the maps disappear. 

Analyzing the results with more detail we can observe in Tables 
\ref{taula contraexemple 1.1} to \ref{taula relacions afins contraexemple} 
that (when the parameter $E$ is small) the map is not 
self-similar, but it behaves ``close to self-similar'' in the 
following sense. 
The values  $\alpha_n'(\omega)$ of the family 
(\ref{FLM C}) (see Tables \ref{taula contraexemple 1.1} to
\ref{taula contraexemple 001.2})
differ form the same values of the family (\ref{FLM B}) 
(see Tables \ref{taula A}  and \ref{taula 2}) an order of magnitude similar to the 
magnitude of $E$. The same happens with the sequences 
$\alpha_n'(\omega)/\alpha_{n-1}'(\omega)$ and  
$\alpha_n'(\omega)/\alpha_{n-1}'(2\omega)$. Actually, this 
indicates that if two maps are ``close'', they still being ``close''
after several renormalizations (or after several magnifications 
of the parameter space).

\section{Summary, conclusions and further development} 

In this paper we have done a numerical study of the asymptotic 
behavior of the slopes $\alpha_n'(\omega,f)$  of the 
reducibility loss bifurcations of quasi-periodic perturbations 
of the Logistic Map. Concretely we have considered families of maps in the 
cylinder $\T\times \R$ which can be written as, 
\begin{equation}
\label{FLM family}
\left.
\begin{array}{rcl}
\bar{\theta} & = & \theta + \omega  ,\\
\bar{x} & = &  \alpha x(1-x) + \eps g(\theta,x)  ,
\end{array}
\right\}
\end{equation}
with $\omega$ an irrational number. Our numerical 
discoveries can be summarized as follows.

\begin{itemize}

\item {\bf First numerical observation} (Section 
\ref{sec:Universality ratio sequence}):  the sequence $\alpha_n'(\omega)/ 
\alpha_{n-1}'(\omega)$ is not convergent in $n$. But, for $\omega$ fix,
one obtains the same sequence for any family of quasi-periodic forced maps 
like (\ref{FLM family}), with
a quasi-periodic forcing of the type $g(\theta,x) =  h(x) \cos(\theta)$.
\item {\bf Second numerical observation} (Section 
\ref{sec:Self renor bif diag}):  the sequence $\alpha_n'(\omega)/ 
\alpha_{n-1}'(2\omega)$ is convergent in $n$ when
the quasi-periodic forcing of the type $g(\theta,x) = h(x) \cos(\theta)$. 
The limit depends on $\omega$
and on the particular family considered.
\item {\bf Third numerical observation} (Section 
\ref{subsection counter example of universality.}): the two previous observations
are not true when the quasi-periodic forcing is of the
type $g(\theta,x) = g_E(\theta,x) = h_1(x) \cos(\theta) + E h_2(x) \cos(2\theta)$
when $E\neq 0$. But the sequence $\alpha_n'(\omega)/ 
\alpha_{n-1}'(2\omega)$ associated to the map (\ref{FLM family})
with $g= g_E$ is $E$-close to the same maps with $g= g_0 $
\end{itemize}

These numerical observations evidence the existence of 
some structure which  govern the asymptotic behavior of the 
sequence $\alpha_n'(\omega)$. This structure is indeed 
the fixed point of a suitable renormalization operator acting 
on the space of functions where the families live. 
The dynamics of this operator determine the asymptotic 
behavior of the sequences $\alpha_n'(\omega)/ 
\alpha_{n-1}'(\omega)$. This dynamics depend on the 
Fourier expansion of the quasi-periodic forcing $g(\theta,x)$, 
giving place to different behaviors depending on the number 
of non-trivial Fourier nodes of $g(\theta,x)$. 
This is described with much more detail in 
the series of papers \cite{JRT11a, JRT11b, JRT11c}. 
In the first one we give the definition of the operator for
the case of quasi-periodic maps and we use it to prove the
existence of reducibility loss bifurcations when the coupling parameter goes to
zero. In the second one we give a theoretical explanation 
to each of the numerical observations above in terms of
the dynamics of the quasi-periodic renormalization operator.
Our quasi-periodic extension of the renormalization
operator is not complete in the sense that
several conjectures must be assumed. In \cite{JRT11c} we include 
numerical computations which support our conjectures and we show
that the theoretical results agree completely with the behavior 
observed numerically.

\bibliographystyle{plain}

\begin{thebibliography}{99}
\bibitem{fn1} 
In \cite{JRT11a} we address this 
problem and we prove that this is the case under suitable hypotheses.
\bibitem{fn2}
Actually we should say $\rho$-geometrically 
asymptotically equivalent, but we will simply say asymptotically equivalent.
\end{thebibliography}

\end{document}